\documentclass{amsart}
\usepackage{amssymb,amsfonts,amsmath}
\usepackage{tikz,youngtab}
\usepackage{hyperref}

\usepackage{color}
\usepackage{filecontents}

\newtheorem{theorem}{Theorem}[section]

\newtheorem{conjecture}[theorem]{Conjecture}

\theoremstyle{definition}

\newtheorem{example}[theorem]{Example}

\numberwithin{equation}{section}
\numberwithin{figure}{section}

\newcommand{\Ehr}{E}   
\newcommand{\Fund}{F}  
\newcommand{\Fstan}{G} 
\newcommand{\ch}{\mathrm{ch}}
\newcommand{\given}{\, | \,}

\definecolor{darkgreen}{rgb}{0,0.5,0}

\title[Stanley's contributions]{The contributions of Stanley to the fabric of 
  \\ symmetric and quasisymmetric functions}

\author{Sara C.  Billey}%
\address{Department of Mathematics, University of Washington, 
Seattle, WA 98195, USA} 
\email{billey@math.washington.edu}%
\urladdr{http://www.math.washington.edu/~billey/}
\thanks{}

\author{Peter R. W. McNamara} \address{Department of Mathematics,
  Bucknell University, Lewisburg, PA 17837, USA}
\email{peter.mcnamara@bucknell.edu}
\urladdr{www.facstaff.bucknell.edu/pm040} 

\thanks{This work was partially supported by a grant from the Simons
  Foundation (\#245597 to Peter McNamara) and by a grant from the National Science Foundation 
  (DMS-1101017 to Sara Billey).}

\date{\today}



\dedicatory{\vspace{1ex}
Dedicated to \ 
$\Yvcentermath1\Yboxdim3ex \  \young(AES,LNT,Y) \ \ \ \young(ADI,CHR,R)$ 
\ on the occasion of his 70th birthday.
\bigskip
}

\begin{document}
{}

\vspace{1in}

\begin{abstract}
We weave together a tale of two rings, SYM and QSYM, following one gold
thread spun by Richard Stanley. The lesson we learn from this tale is
that ``Combinatorial objects like to be counted by quasisymmetric
functions.''
\end{abstract}

\maketitle



\section{Introduction}\label{sec:intro}

The twentieth century was a remarkable era for the theory of symmetric
functions.  Schur expanded the range of applications far beyond roots
of polynomials to the representation theory of $GL_n$ and $S_n$ and
beyond.  Specht, Hall and Macdonald unified the algebraic theory
making it far more accessible.  Lesieur recognized the connection
between Schur functions and the topology of Grassmannian manifolds
spurring interest and further developments by Borel, Bott,
Bernstein--Gelfand--Gelfand, Demazure and many others.  Now, symmetric
functions routinely appear in many aspects of mathematics and 
theoretical physics, and have significant importance in quantum
computation.

In that era of mathematical giants, Richard Stanley's contributions to
symmetric functions are shining examples of how enumerative
combinatorics has inspired and influenced some of the best work of the
century. In this article, we focus on a few of the gems that continue
to grow in importance over time.  Specifically, we survey some results
and applications for Stanley symmetric functions, chromatic symmetric
functions, $P$-partitions, generalized Robinson--Schensted--Knuth
correspondence, and flag symmetry of posets.

As the twentieth century was the century of symmetric functions, then
perhaps the twenty-first century will be defined by the explosion of
developments in the theory of quasisymmetric functions.  The ring of
quasisymmetric functions (QSYM) contains the symmetric functions
(SYM).  It is defined as the subring of power series on countably many
variables with the property that their coefficients are invariant
under arbitrary shifts in the indexing of the variables.  Richard
Stanley first started using quasisymmetric functions in his thesis
while enumerating plane partitions.  Inspired by these results, Ira
Gessel, Richard's first Ph.D. student, recognized the impact of this new tool and formalized the study
of QSYM using the monomial and fundamental bases, relating them to
Schur functions, the Kronecker coefficients and internal products of
skew Schur functions.  Today the quasisymmetric functions are rapidly
growing in importance.  They appear in enumerative combinatorics,
representation theory of $S_{n}$ and 0-Hecke algebras, Macdonald
polynomials and the geometry of Hilbert schemes of points in the
plane, and the Euler--Zagier sums in number theory.

Richard Stanley's contributions to quasisymmetric functions are far
greater than simply giving birth to the field.  He has nurtured the
subject and taught others to expand their influence.  He has at least
36 publications in which symmetric or quasisymmetric functions have
played a role, so we have not attempted to be comprehensive. The main
theme of this paper is to capture Richard's secret to success in this
area.  It starts with simple enumerative questions such as, ``how many
colorings are there for a graph?'', ``how many plane partitions are
there?'', ``how many reduced words exist for a given permutation?'',
or ``how many chains or multichains does a poset have?''.  He relates
the counting problem to a family of generating functions which are
naturally quasisymmetric.  If the quasisymmetric functions are
actually symmetric or even Schur positive, then a beautiful theory
must be at play.  Richard harnesses this symmetric function point of
view to make further advances.

We will assume some familiarity with symmetric and quasisymmetric
functions.  In particular, we follow the notation and terminology of
\cite{Mac95, ec2} to the greatest extent possible.  There is the problem
that $F$ has been overused for quasisymmetric functions.  We will use
$\Fund_{S}(X)$ for the fundamental quasisymmetric function with jumps
in the set $S$.  We rename Ehrenborg's flag generating function for
posets by $\Ehr_{P}(X)$, and we use $\Fstan_{w}(X)$ for Stanley symmetric
functions.

\section{The beginnings}\label{sec:beginnings}

Richard Stanley's work on symmetric functions begins in his Ph.D.\
thesis \cite{StaThesis71}, a subset \cite{StaThesis72} of which is
published in the Memoirs of the AMS series.  Although the discussion
in \cite{StaThesis72} that is directly relevant to symmetric functions
is only one page long, it sets the stage for Ira Gessel's introduction
of quasisymmetric symmetric functions in \cite{Ges84}.  The area of
quasisymmetric functions has seen steady growth since, with particular
acceleration in the last decade.  In Subsection~\ref{sub:quasi}, we explain the
content of this one page of \cite{StaThesis72} to exhibit the role it
played as a catalyst for subsequent developments.

Also while a graduate student,\footnote{Richard recalls writing \cite{Sta71} as a graduate student, some time before the summer of 1970.  See \cite{Sta14} for more details on the timing of Richard's graduation.} Richard wrote \cite{Sta71}, his first paper that had a focus on symmetric functions.  We will elaborate in Subsection~\ref{sub:planepartitions}, explaining how this pair of papers gave an early indication of the flavor of much of Richard's later work on symmetric functions.

\subsection{Quasisymmetric functions}\label{sub:quasi}

Although the term ``quasisymmetric'' was not defined until the eighties \cite{Ges84}, Richard in \cite{StaThesis72} introduces a generating function for $(P,\omega)$-partitions that is a quasisymmetric function.  He also conjectures a condition on $(P,\omega)$ for the generating function to be a symmetric function.  Here, we give the minimal background to explain this generating function, and refer the reader to Gessel's survey in this volume \cite{gessel.tba} for more details. 

For a poset $P$ with $p$ elements, a given bijection $\omega \colon P \to
\{1,2,\ldots,p\}$ can be considered a labeling of the elements of $P$.
If $x<y$ is a covering relation in $P$ and $\omega(x) > \omega(y)$, then we will
call the corresponding edge in the Hasse diagram a $\emph{strict edge}$, otherwise it is a $\emph{weak edge}$.
See Figure~\ref{fig:skewschur}(b) for one example of a labeled poset, where double edges denote strict edges.
A \emph{$(P,\omega)$-partition} $\sigma$ is an order-preserving map
from $P$ to the positive integers that is strictly order-preserving
along strict edges.\footnote{Richard's definition of $(P,\omega)$-partitions in
    \cite{StaThesis72} differs from the one given here in two ways.
    First, Richard gives the nonnegative integers as the codomain, but
    the section of \cite{StaThesis72} of interest to us has the
    positive integers as the codomain. Secondly, his original
    definition has order-reversing in place of order-preserving.  We
    adopt the now customary definition given here since it is what
    Richard uses in \cite{rstan.1995} and our discussion of that paper
    in Section 4 is cleaner if we use the same convention.}
Equivalently, $\sigma$ satisfies the properties:
\begin{enumerate}
\renewcommand{\theenumi}{\alph{enumi}}
\item if $x < y$ in $P$, then $\sigma(x) \leq \sigma(y)$;
\item if $x < y$ in $P$ and $\omega(x) > \omega(y)$, then $\sigma(x) < \sigma(y)$.  
\end{enumerate}

Note that if $P$ is a chain with all weak edges, then a $(P,\omega)$-partition simply corresponds to a partition of a positive integer.  The case when $P$ is a chain with all strict edges gives rise to partitions with distinct parts.  Thus $(P,\omega)$-partitions generalize these classical ideas, hence their name.  Since $\omega$ is a bijection, we can refer to elements of the poset in terms of their $\omega$-labels, and it will be convenient to do so from this point on.  

On page 81 of \cite{StaThesis72}, Richard introduces the generating function in the infinite set of variables $ X =\{x_1, x_2, \ldots\}$ given by
\begin{equation}\label{eq:ppartition}
\Gamma(P,\omega) = \sum_\sigma x_{\sigma(1)} x_{\sigma(2)} \cdots x_{\sigma(p)}\,,
\end{equation}
where the sum is over all $(P,\omega)$-partitions $\sigma$.  This definition is motivated by the fact that skew Schur functions arise as a particular instance of $\Gamma(P,\omega)$, as we explain in the following example.  

\begin{example}
Given a skew diagram $\lambda/\mu$ in English notation with $p$ cells,
label the cells with the numbers $\{1,2,\ldots,p\}$ in any way that
makes the labels increase up columns and from left to right along
rows, as in Figure~\ref{fig:skewschur}(a).  Rotating the result
135$^\circ$ is a counterclockwise direction and replacing the cells
  by nodes as in Figure~\ref{fig:skewschur}(b), we get a corresponding
  labeled poset which we denote by $(P_{\lambda/\mu}, \omega)$ and
  call a \emph{skew-diagram labeled poset}.  Under this construction,
  we see that a $(P_{\lambda/\mu}, \omega)$-partition corresponds
  exactly to a semistandard Young tableau of shape
  $\lambda/\mu$. Therefore $\Gamma(P_{\lambda/\mu},\omega)$ is exactly
  the skew Schur function $s_{\lambda/\mu}$, and is hence a symmetric
  function.  This latter observation appears as Proposition~21.1 in
  \cite{StaThesis72}, at which point Richard states that when $\mu$ is
  empty, $\Gamma(P,\omega)$ is known as a Schur function.
\begin{figure}[htbp]
\begin{center}
\begin{tikzpicture}[scale=0.7]
\begin{scope}[yshift=-0.4cm]
\draw[thick] (0,0) -- (3,0) -- (3,1) -- (4,1) -- (4,3) -- (2,3) -- (2,2) -- (1,2) -- (1,1) -- (0,1) -- cycle;
\draw (1,0) -- (1,1) -- (3,1) -- (3,3)
(2,0) -- (2,2) -- (4,2);
\draw[dashed] (0,1) -- (0,3) -- (2,3)
(0,2) -- (1,2) -- (1,3);
\begin{scope}[font = \Large]
\draw (0.5,0.5) node {1};
\draw (1.5,0.5) node {2};
\draw (2.5,0.5) node {3};
\draw (1.5,1.5) node {4};
\draw (2.5,1.5) node {5};
\draw (3.5,1.5) node {7};
\draw (2.5,2.5) node {6};
\draw (3.5,2.5) node {8};
\end{scope}
\end{scope}
\draw (2,-1) node {(a)};
\begin{scope}[xshift=7cm]
\begin{scope} 
\tikzstyle{every node}=[draw, shape=circle, inner sep=2pt]; 
\draw (5,0) node (a1) {1};
\draw (4,1) node (a2) {2};
\draw (3,2) node (a3) {3};
\draw (3,0) node (a4) {4};
\draw (2,1) node (a5) {5};
\draw (1,0) node (a6) {6};
\draw (1,2) node (a7) {7};
\draw (0,1) node (a8) {8};
\draw[double distance=2pt] (a2) -- (a4)
(a3) -- (a5) -- (a6)
(a7) -- (a8);
\draw (a6) -- (a8)
(a4) -- (a5) -- (a7)
(a1) --(a2) -- (a3); 
\end{scope}
\draw (2.5,-1) node {(b)};
\end{scope}
\end{tikzpicture}
\caption{The skew diagram $443/21$ and a corresponding labeled poset.  Double edges denote strict edges.}
\label{fig:skewschur}
\end{center}
\end{figure}
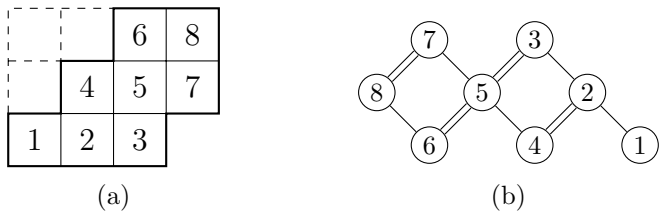
\end{example}

In general, $\Gamma(P,\omega)$ is a quasisymmetric function, meaning that for every sequence of positive integers $a_1, a_2, \ldots, a_k$, the coefficient of $x_{i_1}^{a_1} x_{i_2}^{a_2} \cdots x_{i_k}^{a_k}$ equals that of $x_{j_1}^{a_1} x_{j_2}^{a_2} \cdots x_{j_k}^{a_k}$ whenever $i_1 < i_2 < \cdots < i_k$ and $j_1 < j_2 < \cdots < j_k$.  Note that we get the definition of a symmetric function if we instead allow the $i$'s and $j$'s to be arbitrary sequences of distinct positive integers.  

This leads to what can be called ``Stanley's $P$-partitions conjecture."  An \emph{isomorphism of labeled posets} is a bijection that preserves both the order relation on the poset elements and the set of strict edges.

\begin{conjecture}[\cite{StaThesis72}]
Every finite labeled poset $(P,\omega)$ for which $\Gamma(P,\omega)$ is symmetric is isomorphic to a skew-diagram labeled poset.  
\end{conjecture}

This conjecture remains open; for further information on its status, see \cite{McN06} and the references therein, particularly the Ph.D. thesis work of Malvenuto \cite{Mal93,MalThesis} motivated by Richard's conjecture.

The story resumes in \cite{Ges84}, where Gessel introduces the term
``quasisymmetric,'' with $\Gamma(P,\omega)$ as his motivating example.
Among other things, Gessel defines the two classical bases for the
ring of quasisymmetric functions.  The first is the most natural
  basis of QSYM given by the \emph{monomial quasisymmetric functions}
  $M_{\mathbf{a}}$ indexed by compositions ${\mathbf{a}}=(a_1, a_2, \ldots, a_k)$ where
  $M_{\mathbf{a}}$ is the sum over all monomials of the form $x_{i_1}^{a_1}
  x_{i_2}^{a_2} \cdots x_{i_k}^{a_k}$ with $i_1 < i_2 < \cdots < i_k$.
  The second basis has turned out to be the more important of these
  bases, foreshadowed by its name.   The \emph{fundamental quasisymmetric
    functions} of degree $n$, denoted  $\Fund_S^n(X)$ or just $\Fund_S$ if $n$ and $X$ are understood,  are defined   for $S \subseteq [n-1]$ as follows:
\begin{equation}\label{eq:fund}
\Fund_S^n(X) = \sum_{\genfrac{}{}{0pt}{}{1 \leq i_1 \leq i_2 \leq \cdots \leq i_n}{i_j < i_{j+1} \mathrm{\ if\ } j\in S}} x_{i_1} x_{i_2} \cdots x_{i_n}.
\end{equation}
The well-known bijection from compositions of $n$ to subsets of $n-1$
confirms these two bases are equinumerous in each degree.  
Comparing \eqref{eq:ppartition} and \eqref{eq:fund}, we see the correspondence $\Gamma(P,\omega) =\Fund_S^n$ that occurs when $P$ is a chain of length $n-1$, with $S$ containing $j \in [n-1]$ if and only if the $j$th lowest edge of $P$ is strict.  

Gessel translates Richard's fundamental theorem of
$(P,\omega)$-partitions to the quasisymmetric setting, thus showing the beautiful
and simple expression of $\Gamma(P,\omega)$ for general $P$ in terms of the descent
sets of the linear extensions of $(P,\omega)$:
\begin{equation}\label{eq:gammaf}
\Gamma(P,\omega) = \sum_{\pi \in \mathcal{L}(P,\omega)} F_{D(\pi)},
\end{equation}
where the sum is over all linear extensions $\pi$ of $(P,\omega)$, and $D(\pi)$ denotes the descent set of $\pi$. 

For example, in the labeled poset $(P_{\lambda/\mu}, \omega)$ of Figure~\ref{fig:skewschur}(b), linear extensions are in bijection with SYT of the skew shape in Figure~\ref{fig:skewschur}(a): the SYT record the order in which to take the elements of $(P_{\lambda/\mu}, \omega)$.  We see that we get a descent in the linear extension any time $i+1$ is in a strictly lower row than $i$ in the SYT.  Thus the skew Schur functions expand nicely into the fundamental basis as
\[
s_{\lambda/\mu} = \sum_{T \in \mathrm{SYT}(\lambda/\mu)} F_{D(T)}
\]
where the descent set of a standard Young tableau $T$ of shape $\lambda/\mu$
is the set of all values $i$ such that $i+1$ occurs in a strictly
lower row than $i$. We will see another application of \eqref{eq:gammaf} in
Section~\ref{sub:acyclic}.  

Since 1984, quasisymmetric functions have been developed to such an extent that their importance in algebraic combinatorics is comparable to that of symmetric functions, and we will explore some highlights in the upcoming sections. For more on advances in the area, we refer the reader to 
\cite{Grinberg-Reiner,LMW13} and the many references therein, and to the more informal \cite{BBS10,Wik14}.

\subsection{Theory and Application of Plane Partitions}\label{sub:planepartitions}

In \cite{ec2}, Richard credits Philip Hall \cite{Hal59} with ``the idea of unifying much of the theory of symmetric functions using linear algebra (scalar product, dual bases, involution, etc.).''  In Part 1 of \cite{Sta71}, Richard gives a fuller exposition of the topic of \cite{Hal59}, filling in details and missing proofs.  Richard begins his section on Schur functions by defining them combinatorially in terms of column-strict plane partitions, which is equivalent to the definition in terms of semistandard Young tableaux that is more customary nowadays.  A particularly important thread through the paper is a proof that this combinatorial definition is equivalent to six other definitions of Schur functions: as a determinant, as an orthonormal basis, and in terms of each of the $m$, $h$, $e$ and $p$ bases.\footnote{Both \cite{Hal59} and \cite{Sta71} use the letters $k$, $a$, $e$ and $s$ for bases in place of the letters $m$, $e$, $s$ and $p$, respectively, used in \cite{Mac95, ec2}.}

Although Part 1 of \cite{Sta71} gives missing proofs, it also serves as the first of several examples of Richard's expository writing about symmetric functions.  Undoubtedly, Richard's best known such work is Chapter 7 of \textit{Enumerative Combinatorics} \cite{ec2}.   Other examples are \cite{Sta83}, which was particularly helpful prior to the advent of Appendix~2 of \cite{ec2}, and each of \cite{Sta00,Sta03,Sta04} includes at least one section about symmetric functions. 

Part 2 of \cite{Sta71} exhibits another of Richard's trademarks: using symmetric functions as a tool to prove results that do not involve symmetric functions in their statements.  
We will explore this theme further in the next section.  
We will not describe the content of Part 2 of \cite{Sta71} here, instead referring  the reader to Krattenthaler's survey in this volume \cite{krattenthaler.tba}.

\section{Stanley symmetric functions and applications}\label{sec:ssf+schubert}

Richard found a remarkable application of quasisymmetric functions in
the enumeration of reduced words for permutations.  A \textit{reduced
  word} for a permutation $w$ is a minimal length sequence of positive
integers $\mathbf{a} =(a_1, a_2, \ldots, a_p)$ corresponding to a
product of simple transpositions $s_{a_1} s_{a_2} \cdots s_{a_p} =w$
expressing $w$ in terms of generators $s_i = (i , i+1)$.  A letter in
a reduced word can repeat but not consecutively.  Therefore, reduced
words have  clearly defined ascent sets, denoted $A(a_1,a_2,\ldots,
a_p)=\{i : a_i<a_{i+1}\}\subseteq [p-1]$.  The \textit{Stanley symmetric
  function} for a permutation $w \in S_n$ is defined to be
$$
\Fstan_w(X) = \sum_{(a_1,a_2,\ldots, a_p) \in R(w)} \Fund_{A(a_1,a_2,\ldots, a_p)}(X)
$$
where $R(w)$ is the set of all reduced words for $w$ all of which have
the same length $p=\mathit{inv}(w)$.  Recall that $\mathit{inv}(w)$ denotes the number of \emph{inversion pairs} $(r,s)$, meaning $r<s$ and $w_r > w_s$ where $w=[w_1,\ldots,w_n] \in S_n$.   While these functions $\Fstan_w(X)$ are clearly homogeneous of degree $p$ and quasisymmetric, Theorem 2.1 in
\cite{Sta84} says that $\Fstan_w(X)$ is a symmetric function for all
permutations $w$.\footnote{Richard actually defined
  $F_w(X)=\sum_{(a_1,a_2,\ldots, a_p) \in R(w)}
  \Fund_{D(a_1,a_2,\ldots, a_p)}(X)$ which equals $G_{w^{-1}}$ in our
  notation.  The switch from $w$ to $w^{-1}$ is related to a formula
  for Schubert polynomials coming in Subsection~\ref{sub:schubert}.}  The
original proof exhibits symmetry via an intricate bijection among the
reduced words contributing to the coefficients of
$x_1^{\alpha_1}\cdots x_k^{\alpha_k}$ and $x_1^{\alpha_1}\cdots
x_{i}^{\alpha_{i+1}} x_{i+1}^{\alpha_{i}}\cdots x_k^{\alpha_k}$.

Richard identified a special family of $\Fstan_w$'s in Section 4 of
\cite{Sta84}.  In particular, $w$ is 2143-avoiding (vexillary) if and
only if $\Fstan_w = s_{\lambda(w)}$ where $\lambda(w)$ is the
partition obtained by sorting the Lehmer code of $w$ which is the
sequence $(c_1,c_2,\ldots, c_n)$ with $c_i=|\{j>i \given w_i>w_j\}|$.
The vexillary permutations were an early application of pattern
avoidance.  Pattern avoidance has grown into an important and very active research
area.  One of the main problems in that area was the Stanley--Wilf
conjecture from around 1980 which says that for any permutation $w$, the number of
permutations in $S_n$ avoiding $w$ is at most $c^n$ for some constant
$c$.  This conjecture was proved in 2004 by Marcus and Tardos
\cite{marcus.tardos.2004}.  See Richard's portion of \cite{CGH15} for more on the genesis and development of the Stanley--Wilf conjecture.  

Richard's initial interest in these functions $\Fstan_w$ was for counting the
number of reduced expressions for any $w$, in particular for the
longest permutation $w_0=[n,n-1,\ldots,1] \in S_n$.  Note $|R(w)|$ is
equal to the coefficient of $x_1 x_2 \cdots x_p$ in $\Fstan_w$
provided $w$ has $p$ inversions.  Because $w_0$ is vexillary, he was
able to show that $|R(w_0)|$ is the number of standard tableaux of
staircase shape $(n-1,n-2,..., 1)$, which is easily computed via the
Frame--Robinson--Thrall hook length formula.  As a function of $n$,
this sequence grows very fast: 1, 1, 2, 16, 768, 292864.  See
\cite[A005118]{oeis}.  The other vexillary permutations $v$ have
similarly easy formulas: $|R(v)|=f^{\lambda(v)}$ where $f^{\lambda}$
counts the number of standard Young tableaux of shape $\lambda$.  More
generally, Richard conjectured that every $\Fstan_w$ was \textit{Schur
  positive}: $G_{w} = \sum a_{\lambda ,w} s_{\lambda}$ with
$a_{\lambda,w} \in \mathbb{N}$.  Thus, the expansion coefficients
$a_{\lambda,w}$ could be used to calculate $|R(w)| = \sum
a_{\lambda,w}f^{\lambda}$.

Edelman--Greene \cite{EG} proved Richard's conjecture shortly
thereafter.  See also Lascoux--Sch\"utzenberger \cite{LS7} on the
plactic monoid for an alternative approach.  The Edelman--Greene
correspondence is a slightly modified version of the classical RSK
algorithm for reduced words: when inserting an $i$ into a row $j$
already containing an $i$, skip row $j$, and enter into the next row
the larger of $i,i+1$ occurring in row $j$.  Each reduced word
$\mathbf{a}$ bijectively gives rise to a pair of tableaux
$(P(\mathbf{a}),Q(\mathbf{a}))$ known as the insertion tableau and the
recording tableau respectively.  Edelman--Greene show that for each
insertion tableau $P$ that arises when doing their insertion algorithm
on all reduced words for $w$ and every standard tableau $Q$ of the
same shape as $P$, there exists a unique $\mathbf{a} \in R(w)$ with
$P(\mathbf{a})=P$ and $Q(\mathbf{a}) =Q$.  Thus, $a_{\lambda,w}$
counts the number of distinct $P$ tableaux of shape $\lambda$ that
arise from $w$.  

There is another beautiful proof that all of the Stanley symmetric
functions are symmetric and Schur positive.    Little gave a bijection termed
a ``bumping algorithm'' on reduced words which preserves ascent sets
thus proving an effective recurrence
\begin{equation}\label{e:trans}
 \Fstan_w = \sum_{w' \in T(w)}  \Fstan_{w'}
\end{equation}
terminating when $w$ is vexillary in which case
$\Fstan_w=s_{\lambda(w)}$ \cite{little2003combinatorial}.  Here $T(w)
:=\{w'=vt_{i,r} \given i<r, \mathit{inv}(w')=\mathit{inv}(w)\}$ where
$t_{a,b}$ is the transposition interchanging $a$ and $b$, and $v=wt_{r,s}$ with $(r,s)$ being the lexicographically largest inversion pair.  The recurrence in
\eqref{e:trans} is known as the \textit{transition equation} for
Stanley symmetric functions and was originally proved by
Lascoux--Sch\"utzenberger \cite{LS.1985} who suggest it can be used to
compute Littlewood--Richardson coefficients effectively.

Recently, Hamaker and Young \cite{HamakerYoung2014} proved a
conjecture of Thomas Lam's that Little bumps preserve the $Q$
tableaux. Thus, the $a_{\lambda,w}$ also count the number of reduced
words for $w$ in the same communication class under Little bumps as a
fixed reduced word for the unique permutation up to trailing fixed
points with Lehmer code $(\lambda_1, \lambda_2, \ldots,
\lambda_k,0,0,\ldots)$.

Stanley symmetric functions and the enumeration of reduced words are
now known to have many applications and connections to representation
theory of $S_n$, geometry of Schubert varieties, and stochastic
processes related to sorting networks.  For example, via the work of
Kra{\'s}kiewicz \cite{kraskiewicz} and Reiner--Shimozono
\cite{plactification}, there is a generalization of a Specht module on
the (Rothe) diagram of the permutation of $w$ which has $\Fstan_w$ as
its Frobenius characteristic.  Also, Pawlowski showed that the cohomology classes of Coskun's
rank varieties in the Grassmannian manifolds are all Stanley symmetric functions \cite{paw.2014}.

\subsection{Random reduced words}

Angel--Gorin--Holroyd--Romik--Vir\'ag \cite{angel2007random,angel2012pattern} have initiated a program to study random
reduced expressions for the longest permutation
 and related processes.  They
produce a random reduced expression uniformly by using the hook walk
algorithm due to Greene--Nijenhuis--Wilf \cite{Greene-Nijenhuis-Wilf} to produce a uniformly random
staircase shape standard tableau $Q$, and apply the inverse
Edelman--Greene correspondence along with the unique $P$ tableau for
this permutation.  The hook walk algorithm uniformly at random chooses
one cell in the shape $\lambda$, then it ``walks'' to a different cell
in the hook of this first cell uniformly.  From that cell, it again
chooses a new cell in its hook to walk to, continuing until the walk
arrives at a corner cell, placing the largest value there. Then restricting to the still empty cells of $\lambda$, the algorithm repeats the process to place the second largest value, etc.    The following tantalizing conjecture is
still open at this time.

\begin{conjecture}\cite{angel2007random}\label{conj:AHRV}
  Choose a uniform random reduced word $(a_1,\ldots, a_p)$ for $w_0
  \in S_n$.  The probability distribution of the 1's in the
  permutation matrix for the initial product $s_{a_1}s_{a_2} \cdots s_{ a_{\lfloor p/2 \rfloor }}$
  approaches the surface measure of the sphere projected
  to 2 dimensions as $n$ gets large.
\end{conjecture}
\begin{figure}
 \includegraphics[height=3cm]{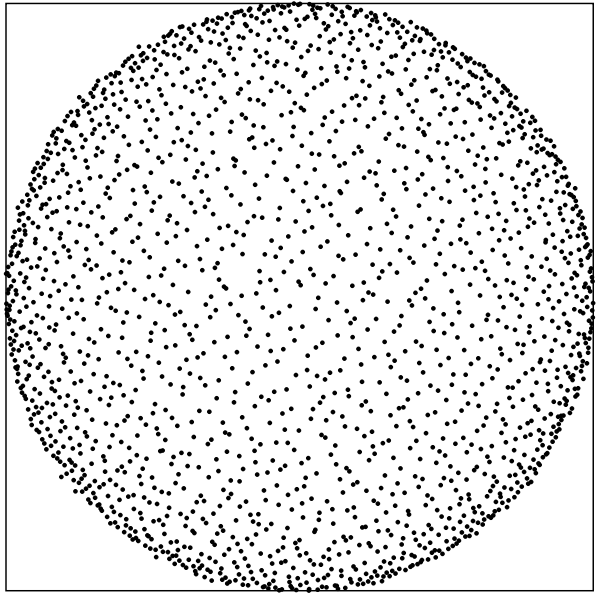} 
  \caption{An example of a randomly chosen permutation matrix with $n=2000$ demonstrating the Angel-Holroyd-Romik-Vir\'ag Conjecture~\ref{conj:AHRV} from \cite{angel2007random}. }
\end{figure}  

\subsection{Schubert polynomials}\label{sub:schubert}

Schubert polynomials $\mathfrak{S}_w$ for $w \in S_n$ are a
generalization of Schur polynomials that were invented by Lascoux and
Sch\"utzenberger in the early 1980s.  Schur polynomials represent the
Schubert basis for the cohomology ring of a Grassmannian manifold.
This connection is realized by comparing the Giambelli formula with
the Jacobi--Trudi determinantal formula \cite{Fulton-book}.  Using the
Bernstein--Gelfand--Gelfand divided difference operators \cite{BGG},
Lascoux and
Sch\"utzenberger showed that Schubert polynomials are explicit representatives of
the cohomology classes of Schubert varieties in the complete flag
manifold $GL_n/B$ where $B$ is the set of upper triangular matrices.
They form a basis for the cohomology ring which can be identified as
polynomials in $n$ variables modulo the ideal of  symmetric polynomials with no constant term.
Schubert polynomials have two distinct advantages over
other such bases. First, their structure coefficients determine the
cup product for the cohomology ring without ever having to reduce
modulo the ideal of symmetric polynomials.  Second, they have a
positive expansion into monomials. Lascoux and Sch\"utzenberger
pointed out that Stanley symmetric functions are the limiting case of
their Schubert polynomials in the sense that $\Fstan_w =
\lim_{k\longrightarrow \infty} \mathfrak{S}_{1^k \times w}$ where
$1^k\times w = [1,2,\ldots, k, w_1+k, \ldots, w_n+k]$.

Around 1991, Richard conjectured the following explicit formula for
Schubert polynomials, formalizing the connection to his symmetric
functions and arguably revolutionizing the field:
\begin{equation}\label{e:bjs}
\mathfrak{S}_w = \sum_{\mathbf{a}=(a_1,\ldots,a_p) \in R(w)} \ 
\sum_{(i_1,\ldots,i_p)\in C(\mathbf{a})} x_{i_1}x_{i_2} \cdots x_{i_p}
\end{equation}
where $C(\mathbf{a})$ are the \textit{compatible sequences} of
$\mathbf{a}$ defined very similarly to the terms in the fundamental
quasisymmetric functions.  Specifically, $(i_1,\ldots,i_p)\in
C(\mathbf{a})$ provided

\begin{enumerate}
\item $i_1  \leq i_2 \leq \ldots \leq i_p$ are positive integers,
\item if $a_j<a_{j+1}$ then $i_j<i_{j+1}$, and 
\item each $i_j \leq a_j$ for all $1\leq j\leq p$.
\end{enumerate}
For example, $w=[1,4,3,2]$ has two reduced words $(3,2,3)$ and
$(2,3,2)$ and they have 5 compatible sequences: 
$C(3,2,3) =\{(1,1,2), (1,1,3), (1,2,3),(2,2,3)\}$ and $C(2,3,2)=\{(1,2,2)\}$.  So
$$
\mathfrak{S}_{[1,4,3,2]} =  x_1^2 x_2 + x_1^2x_3 + x_1x_2x_3 + x_2^2x_3 + x_1x_2^2.
$$

\begin{figure}
\begin{tikzpicture}[scale=0.3]
\begin{scope}
\draw (0,0) -- (4.8,4.8) 
(5.2,5.2) --  (5.8,5.8) --  (6.2,5.8) -- (10,2);
\draw (0,2) -- (3.8,5.8)
(4.2,6.2) -- (4.8,6.8)
(5.2,7.2) -- (5.8,7.8)
(6.2,8.2) -- (6.8,8.8) -- (7.2,8.8) -- (7.8,8.2) -- (8.2,8.2) -- (8.8,8.8)
(9.2,9.2) -- (10,10);
\draw (0,4) -- (2.8,6.8) -- (3.2,6.8) -- (10,0);
\draw (0,6) -- (1.8,7.8) -- (2.2,7.8) -- (2.8,7.2) -- (3.2,7.2) -- (3.8,7.8) -- (4.2,7.8) -- (5.8,6.2) -- (6.2,6.2) -- (6.8,6.8) -- (7.2,6.8) -- (10,4);
\draw (0,8) -- (0.8,8.8)
(1.2,9.2) -- (1.8,9.8) -- (2.2,9.8) -- (2.8,9.2) -- (3.2,9.2) -- (3.8,9.8) -- (4.2,9.8) -- (4.8,9.2) -- (5.2,9.2) -- (5.8,9.8) -- (6.2,9.8) -- (6.8,9.2) -- (7.2,9.2) -- (7.8,9.8) -- (8.2,9.8) -- (10,8);
\draw (0,10) -- (1.8,8.2)-- (2.2,8.2) -- (2.8,8.8) -- (3.2,8.8) -- (3.8,8.2) -- (4.2,8.2) -- (4.8,8.8) -- (5.2,8.8) -- (6.8,7.2) -- (7.2,7.2) -- (7.8,7.8) -- (8.2,7.8) -- (10,6);
\draw (10.5,0) node {\small{3}};
\draw (10.5,2) node {\small{1}};
\draw (10.5,4) node {\small{4}};
\draw (10.5,6) node {\small{6}};
\draw (10.5,8) node {\small{5}};
\draw (10.5,10) node {\small{2}};
\end{scope}
\begin{scope}[xshift=50em, scale=0.9]
\draw (-0.5,8) arc (270:360:0.5) 
(0,8.5) -- (0,11.5);
\draw (-0.5,0) arc (270:360:0.5) 
(0,0.5) -- (0,3.5) 
(0,3.5) arc (180:90:0.5)
(0.5,4) -- (1.5,4)
(1.5,4) arc (270:360:0.5)
(2,4.5) -- (2,11.5);
\draw (-0.5,10) -- (3.5,10)
(3.5, 10) arc (270:360:0.5) 
(4,10.5) -- (4,11.5);
\draw (-0.5,6) arc (270:360:0.5)
(0,6.5) -- (0,7.5)
(0,7.5)  arc (180:90:0.5)
(0.5,8) -- (3.5,8)
(3.5,8) arc (270:360:0.5)
(4,8.5) -- (4,9.5)
(4,9.5) arc (180:90:0.5)
(4.5,10) -- (5.5,10)
(5.5,10) arc (270:360:0.5)
(6,10.5) -- (6,11.5);
\draw (-0.5,2) -- (1.5,2)
(1.5,2) arc (270:360:0.5)
(2,2.5) -- (2,3.5) 
(2,3.5) arc (180:90:0.5)
(2.5,4) -- (3.5,4)
(3.5,4) arc (270:360:0.5)
(4,4.5) -- (4,5.5) 
(4,5.5) arc (180:90:0.5)
(4.5,6) -- (5.5,6)
(5.5,6) arc (270:360:0.5)
(6,6.5) -- (6,7.5) 
(6,7.5) arc (180:90:0.5)
(6.5,8) -- (7.5,8)
(7.5,8) arc (270:360:0.5)
(8,8.5) -- (8,11.5);
\draw (-0.5,4) arc (270:360:0.5)
(0,4.5) -- (0,5.5)
(0,5.5) arc (180:90:0.5)
(0.5,6) -- (3.5,6)
(3.5,6) arc (270:360:0.5)
(4,6.5) -- (4,7.5)
(4,7.5) arc (180:90:0.5)
(4.5,8) -- (5.5,8)
(5.5,8) arc (270:360:0.5)
(6,8.5) -- (6,9.5)
(6,9.5) arc (180:90:0.5)
(6.5,10) -- (9.5,10)
(9.5,10) arc (270:360:0.5)
(10,10.5) -- (10,11.5);
\draw (-1,0) node {\small{2}};
\draw (-1,2) node {\small{5}};
\draw (-1,4) node {\small{6}};
\draw (-1,6) node {\small{4}};
\draw (-1,8) node {\small{1}};
\draw (-1,10) node {\small{3}};
\end{scope}
\end{tikzpicture}
 \caption{A pseudo line arrangement  for $w=[3,1,4,6,5,2]$ is shown on the left and a reduced pipedream for $w$ is shown on the right.  Both encode the RC-pair with $\mathbf{a} =(5,2,1,3,4,5)$ and $\mathbf{i}=(1,1,1,2,3,5)$.}
\label{fig:rcgraphs}
\end{figure}
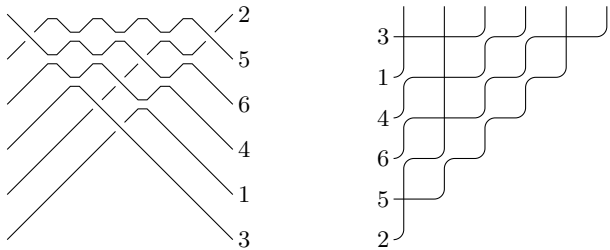

The Schubert polynomial conjecture was proved in two papers in quick
succession \cite{BJS,FS} in 1992.  The pairs of reduced words and
compatible sequences were then restated geometrically in terms of
pseudo line arrangements in \cite{FK} and in terms of RC-graphs in
\cite{BB}, which are also known as \textit{reduced pipedreams} because
of their visual similarity to a game with the same name made in 1989 for the
Commodore Amiga \cite{knutson-miller-2005,wiki:pipedream}.  See Figure~\ref{fig:rcgraphs}.  
 
A reduced pipedream is constructed from a reduced word
$\mathbf{a} = (a_1,\ldots,a_p)$ and one of its compatible sequences
$\mathbf{i} = (i_1,\ldots,i_p)$ by placing a tile with a string crossing in each
entry of a matrix of the form $(i_j, a_j-i_j+1)$ and a tile with two
elbows in all other positions $(i,j)$ with $i+j\leq n$.  
Conversely, any
placement of the tiles taking the strings labeled $1,2,\ldots,n$ along
the top to the labeling $w_1, w_2,\ldots, w_n$ down the left side such
that no pair of strings cross more than once gives rise to a reduced
pipedream for $w$. The corresponding RC-pair $(\mathbf{a}, \mathbf{i})$  
can be recovered from the row and
column numbers of the crossings reading along rows right to left, from
top to bottom.  Thus, the RC-pairs for $w$ or equivalently the
reduced pipedreams $RP(w)$ play the role of semistandard tableaux for
Schur functions. The Schubert polynomial $\mathfrak{S}_w$ is the
generating function
$$
\mathfrak{S}_w(X) = \sum_{D \in RP(w)}  x^D
$$
where $x^D=x_1^{n_1}x_2^{n_2} \cdots$ and for each $i$, $n_i$ is the
number of crossings on row $i$ in $D$. Similarly, the double Schubert
polynomial $\mathfrak{S}_w(X,Y)$ is the generating function
$$
\mathfrak{S}_w(X,Y)= \sum_{D \in RP(w)}  \prod_{(i,j) \in D} (x_i - y_j).
$$
due to a theorem of Fomin--Kirillov \cite{FK}.   
The double Schubert polynomials can be used to represent Schubert
classes in the equivariant cohomology ring of the flag manifold \cite{B4,Graham.2001}.
Recently, this field is moving toward more exotic cohomology theories.
For example, see both \cite{anderson.chen.2014,lam.shimo.2014} for a
quantum equivariant analog.

The individual RC-pairs are also related to the geometry of
Schubert varieties.  In his 2000 Ph.D.\ thesis \cite{Kogan.phd}, Kogan
gave a degeneration of a Schubert variety to a toric variety by
interpreting each of the RC-pairs as inequalities on faces in a
polytope.

In 2005, Knutson and Miller gave a second geometric interpretation to
the reduced pipedreams (or equivalently RC-pairs) using matrix
Schubert varieties $X_w$ in their paper on ``Gr\"obner geometry of
Schubert polynomials'' \cite[Thm.~B]{knutson-miller-2005}, thereby establishing
the pipedreams as the most natural presentation of the RC-pairs.
Let $I_w$ be the determinantal ideal defining $X_w$ by rank
conditions.  The initial ideal $J_w= \mathrm{in}(I_w)$ under an ``antidiagonal
term order'' such as $z_{1,1}>z_{2,1} > \cdots > z_{n,1}> \cdots >
z_{n,n}$ has square-free generators so it corresponds with a
Stanley--Reisner simplicial complex $L_w$.  They show that the
Schubert class $[X_w]$ is equal to the class $[L_w]$.  The class
$[L_w]$ is particularly easy to compute because it is just a union of
coordinate subspaces of $M_n$ which are in bijection with reduced
pipedreams and the explicit equations are determined by the crossings
in the pipedreams.  Furthermore, the prime decomposition of $J_w$ is
given by intersecting the prime ideals $\langle z_{ij} \ | \ (i,j)
\in D\rangle$ for all $D \in RP(w)$.  This interpretation covers both
the connections to cohomology and K-theory for single and double
Schubert polynomials and Grothendieck polynomials.

The study of reduced words and Schubert varieties naturally extends to
other Lie types.  For the other classical groups of types $B$, $C$,
and $D$, the Stanley symmetric functions have been generalized and
used to give analogs of Schubert polynomials \cite{BH,FK2,Lam1}.
Kra{\'s}kiewicz found the analog of the Edelman--Greene correspondence
for type $B$ \cite{kraskiewicz}, and Richard's student Tao-Kai Lam did
the same for type $D$ \cite{Lam1}.  One of the biggest breakthroughs
in Schubert theory of this century is due to Thomas Lam's definition
of the affine Stanley symmetric functions and his theorem \cite{Lam06} that these
are one form of the $k$-Schur functions defined by
Lascoux--Lapointe--Morse \cite{LLM2003}.

The study of
Stanley symmetric functions is a perfect example of how Richard
influences the field.  He started by asking ``how many reduced words
are there for a permutation?'' and he ended up inspiring so much more.
His 1984 paper has been cited now over 70 times according to
MathSciNet.  It certainly counts as a gem.

\section{Chromatic symmetric functions and quasisymmetric functions}\label{sec:chromatic}

What is a good example of a symmetric function that arises naturally
other than the usual bases for the ring of symmetric functions?  Take
a simple graph $G=(V,E)$ with $n$ vertices labeled $V=\{1,2,\dotsc,n
\}$ and consider all of its \textit{proper colorings} $c: V
\longrightarrow \mathbb{P}$, meaning that no two adjacent vertices get
assigned to the same value in $\mathbb{P}$.  Given any proper coloring
of $G$, if we permute the ``colors'' in $\mathbb{P}$ we get another
proper coloring.  The symmetry inherent in graph coloring can be
formalized by the \textit{chromatic symmetric function} for $G$
defined by
\begin{equation}\label{eq:chromatic}
X_{G}(x_{1},x_{2},\ldots ) = \sum_{c:V\longrightarrow \mathbb{P}}
x_{c(1)} x_{c(2)}\cdots x_{c(n)}.
\end{equation}

This colorful family of symmetric functions and some of its
generalizations were invented and first studied by Richard in the
1990s in two papers \cite{rstan.1995,rstan.1998}.  Much of this
work was inspired by properties of the famous chromatic polynomials
$\chi_{G} (m)$ invented by George David Birkhoff and extended by
Hassler Whitney in the early 1900s.  In fact, chromatic polynomials
have been inspiring Richard throughout his career
\cite{rstan.1970,rstan.1973,rstan.2012}.  We will review some of these
connections between $\chi_G$ and $X_G$.  A refinement of Richard's
chromatic symmetric function is a central topic of the paper of Shareshian and
Wachs in this volume \cite{shareshian.wachs.tba}.

Recall that the \textit{chromatic polynomial} $\chi_{G} (m)$ counts
the number of proper colorings of $G$ using $m$ colors.  The proof
that this function is a polynomial in $m$ proceeds simply by
considering a recurrence using the deletion and contraction of any one
edge, or by realizing that $\chi_G(m) = \sum_{i=0}^{|V(G)|}
\binom{m}{i} C_{G}(i)$, where $C_{G}(i)$ is the number of ways to properly color $G$ using \emph{exactly} $i$ colors. We get a third proof of this polynomiality property
by evaluating the chromatic symmetric function at $1^{m}$, which means we set 
$x_{1}=x_{2}=\cdots = x_{m}=1$ and $x_{i}=0$ for all $i>m$.  Observe that $X_{G}(1^{m}) =
\chi_{G}(m)$, and every symmetric function evaluated at $1^{m}$ is a
polynomial function of $m$, as Chow points out in \cite{chow.1999}, which 
is easily verified on generators such as $p_{i}(1^{m})=m$.

Alternating sums for chromatic polynomials and chromatic symmetric functions are
a good example of the interplay between the two subjects.  Whitney's
expansion of the chromatic polynomial is the alternating sum
\begin{equation}
\chi_{G}(m)= \sum_{S \subseteq E(G)} (-1)^{|S|} \ m^{c(S)}
\end{equation}
over all spanning subgraphs of $G$ with edge sets $S$ where $c(S)$
counts the number of connected components of the subgraph
\cite{whitney.1932}.  Richard proved the analogous statement for the
chromatic symmetric functions giving their expansion into power sum
symmetric functions:
\begin{equation}
X_{G}= \sum_{S \subseteq E(G)} (-1)^{|S|} \ p_{\lambda (S)}
\end{equation}
where the sum is again over all spanning subgraphs of $G$ and $\lambda
(S)$ is the partition of $n=|V(G)|$ determined by the sizes of the
connected components of the spanning subgraph with edge set $S$
\cite[Thm 2.5]{rstan.1995}.

\subsection{Acyclic orientations}\label{sub:acyclic}

Acyclic orientations are another good example of how one can extend 
results for $\chi_G$ to the $X_G$ setting.  Interestingly,
Richard's proofs use the $(P,\omega)$-partition ideas of
Subsection~\ref{sub:quasi}. Given a simple graph $G$, choose a
direction for each edge.  If there is no directed cycle in the chosen
orientation on the edges, we say the orientation is
\textit{acyclic}. Let $a_G$ be the number of acyclic orientations of
$G$.  Richard had shown in previous work \cite{rstan.1973} that
\begin{equation}\label{e:acyclic}
  a_G = (-1)^{n}\chi_{G}(-1).
\end{equation}
  This result is often called ``a classic'' which has many
  consequences itself.  For example, Hanlon \cite{hanlon.2008} gave a
  topological explanation for this result using the Hodge
  decomposition of the coloring complex due to Steingr\'imsson
  \cite{steingrimsson.2001}.  See Propp's paper in this volume
  \cite{propp.tba} for further discussion of the enumerative
  consequences.  Observe that every proper coloring gives rise to a
  specific acyclic orientation $\mathfrak{o}$ of the edges of $G$ by orienting from
  the larger colored vertex to the smaller.  Each such acyclic orientation
   gives rise to a poset $P_{\mathfrak{o}}$ on the
  vertices of $G$ by taking the transitive closure of this directed acyclic 
  graph, where orientations point downwards in the Hasse diagram.
Endow $P_{\mathfrak{o}}$ with a labeling $\omega_{\mathfrak{o}}$ that makes all edges strict, as we defined in Subsection~\ref{sub:quasi}.  Relating the definitions \eqref{eq:ppartition} and \eqref{eq:chromatic} and applying \eqref{eq:gammaf}, Richard \cite[p.~176]{rstan.1995} deduced the following expansion of the chromatic symmetric function in terms of
  fundamental quasisymmetric functions:
 \begin{equation}\label{e:chrom.fundamentals}
    X_G = \sum_{\mathfrak{o}} \Gamma(P_{\mathfrak{o}}, \omega_{\mathfrak{o}}) = \sum_{\mathfrak{o}}\ \sum_{\pi \in \mathcal{L}(P_{\mathfrak{o}}, \omega_{\mathfrak{o}})}  \Fund_{D(\pi)}
  \end{equation}
where the first sum is over all acyclic orientations and the second
sum is over all linear extensions of the corresponding labeled posets.  

Chow \cite{chow.1999} shows how
\eqref{e:chrom.fundamentals} implies the following result of Chung and
Graham \cite[Thm.~2]{chung.graham.1995}.  If a chromatic polynomial
$\chi_{G}(x) $ is expanded in the basis of binomial coefficients  
\[
\binom{x+k}{n}  \text{ for } k=0,1,\ldots 
\]
then the coefficients are nonnegative and count the number of
bijective labelings $\pi$ of the vertices of $G$ with exactly $k$ of
what they call $G$-descents.

Richard uses \eqref{e:chrom.fundamentals} to give what he calls the
``main result on acyclic orientations'' connecting the expansion of
$X_G$ in the elementary symmetric function basis with acyclic
orientations.  Say
\begin{equation}\label{e:chromatic.elementary} X_G =
\sum_{\lambda \vdash n} c_{\lambda} e_{\lambda},
\end{equation}
and let $\mathrm{sink}(G,j)$ be the number of acyclic orientations of
$G$ with exactly $j$ sinks.  Then \cite[Thm.~3.3]{rstan.1995} tells
us
\begin{equation}\label{e:sinks}
\mathrm{sink}(G,j) = \sum_{\substack{\lambda \vdash n \\l(\lambda ) =j
}} c_{\lambda },
\end{equation}
where $l(\lambda)$ denotes the number of parts of $\lambda$.
The proof is via a linear transformation defined on fundamental
quasisymmetric functions by mapping $\Fund_{S}$ to $t(t-1)^{i}$ if
$S=\{i+1,i+2,\ldots, n-1 \}$ and 0 otherwise.  Richard uses \eqref{e:chrom.fundamentals}
to show that the transformation applied to the left-hand side of \eqref{e:chromatic.elementary} yields
$\sum_j \mathrm{sink}(G,j)t^{j}$.  He again uses the idea of $(P,\omega)$-partitions to show that the transformation applied to the right-hand side of \eqref{e:chromatic.elementary} yields $\sum_\lambda c_{\lambda} t^{l(\lambda )}$ where the sum is over the desired $\lambda$.

\subsection{Combinatorial Hopf algebras}

Among the most important theorems in quasisymmetric function theory
is that the ring of quasisymmetric functions is the terminal object
in the category of combinatorial Hopf algebras.  This foundational
result is due to Aguiar--Bergeron--Sottile \cite{ABS}.  They claim it
explains the ubiquity of quasisymmetric functions in mathematics.

A \textit{combinatorial Hopf algebra} is a graded connected Hopf
algebra $\mathcal{H}$ over a field $\mathbb{F}$ along with a choice of
character $\zeta : \mathcal{H} \longrightarrow \mathbb{F}$.  The ring
of quasisymmetric functions is a combinatorial Hopf algebra with the
canonical character $\zeta_Q$ defined on the monomial basis by saying 
$\zeta_Q(M_\alpha)$ is 1 or 0 depending on whether or not $\alpha$ is
a composition with at most 1 part.   

A classic example is Schmitt's Hopf algebra on finite graphs
\cite{Sch94}.  Let $\mathcal{H_G}$ be the $\mathbb{F}$-vector space
with basis given by the isomorphism classes of finite simple graphs.
Given a graph $G$ and a subset of the vertices $S \subseteq V(G)$, let
$G|_{S}$ be the induced subgraph on the vertices in $S$.
Multiplication of graphs in $\mathcal{H_G}$ is given by disjoint union
and comultiplication is given by
\[
\Delta(G) = \sum_{S {\subseteq} V (G)} G|_{S} \otimes G|_{V(G) \setminus S}.
\]
A character on $\mathcal{H_G}$ can be defined by $\zeta (G) =1$ if $E(G)
=\emptyset$ and 0 otherwise.  Note that a graph with no edges can have
its vertices all colored the same in a proper coloring.
The Aguiar--Bergeron--Sottile theorem specifies the  explicit 
morphism of combinatorial Hopf algebras from $(\mathcal{H_G}, \zeta)$
to $(\mathrm{QSYM}, \zeta_Q)$.  It maps any basis element $G \in
\mathcal{H_G}$ to $\sum \zeta_{\mathbf{a}} (G) M_{\mathbf{a} }$ where
the sum is over all compositions $\mathbf{a}$ of $n=|V(G)|$ and $M_{\mathbf{a} }$
is the monomial quasisymmetric function.  Since $\zeta_{\mathbf{a}}
(G)$ counts the number of ways to partition $G$ into edgeless spanning
subgraphs of sizes $\mathbf{a} = (a_{1}, a_{2},\dotsc , a_{k})$, we
see that $G$ in Schmitt's Hopf algebra maps to the chromatic symmetric function $X_G$.  

\subsection{Open problems} 

In typical Stanley style, we close this section with some of the many
interesting open problems related to colorings of graphs.

\begin{enumerate}
\item Which polynomials with integer coefficients are
chromatic polynomials\linebreak \cite{Rea68}?

\item When do two graphs have the same chromatic polynomial?  

\item When do two graphs have the same chromatic symmetric function? 

\item If two trees have the same chromatic symmetric function, are
  they necessarily isomorphic \cite{rstan.1995}?  Perhaps so.  At
  least it holds for trees with up to 23 vertices as checked by
  Li-Yang Tan.  See also Martin--Morin--Wagner
  \cite{Martin-Morin-Wagner} and the references therein for further partial results, and \cite{AZ14, OS14, SST15} for more recent progress.

\item Are the chromatic symmetric functions of incomparability graphs
  of $(3+1)$-free posets $e$-positive?  Richard and John Stembridge conjecture yes
  \cite[Conj.~5.5]{SS93}\cite[Conj.~5.1]{rstan.1995} and give supporting evidence.
  Gasharov proved that in this case, $X_G$ is Schur positive
  \cite{Gas96}.  Recently, there has been some exciting progress on
  this conjecture due to Mathieu Guay-Paquet.  He has reduced the
  problem to the subclass of $(3+1)$-and-$(2+2)$-free posets, which
  are called ``semiorders'' or ``unit interval orders'' in the
  literature and are enumerated by Catalan numbers \cite{GP.2013}.
  See also the alternative approach and further conjectures by
  Shareshain and Wachs using representation theory on the cohomology
  groups of Hessenberg varieties \cite{2014.shareshian.wachs,shareshian.wachs.tba}.  
\end{enumerate}

\section{A skew generalization of the RSK algorithm}\label{sec:rsk}

A fundamental result in the theory of symmetric functions is the Robinson--Schensted--Knuth (RSK) algorithm \cite{Rob38,Sch61,Knu70}, which gives a bijection between matrices $A$ over $\mathbb{N}$ with a finite number of non-zero entries, and pairs $(P,Q)$ of semistandard Young tableaux (SSYT) of the same shape.  Among the consequences of the RSK algorithm (see \cite[\S 7.12]{ec2} for more) is the Cauchy identity for Schur functions, from which the orthonormality of the basis of Schur functions follows.  Restricting $A$ to the case of $n\times n$ permutation matrices, the bijectivity of the RSK algorithm implies the beautiful identity
\begin{equation}\label{equ:flambda}
\sum_{\lambda \vdash n} (f^\lambda)^2 = n!\,,
\end{equation}
where $f^\lambda$ denotes the number of standard Young tableaux (SYT) of shape $\lambda$.

As one would expect, there has been much work done in developing analogues and generalizations of the RSK algorithm; a brief overview of such work can be found in the introduction to \cite{SS90}.  Our goal for this section is to highlight the generalization of the RSK algorithm to skew shapes due to Bruce Sagan and Richard \cite{SS90}, and point out some more recent applications of their generalization.   Another perspective on \cite{SS90} and its connection to representation theory appears in Lenart's paper in this volume \cite{lenart.tba}.

\subsection{The skew version of the RSK algorithm}  

The classical RSK algorithm works by first converting the matrix $A=(a_{ij})$ to a word that consists of $a_{ij}$ copies of the biletter $(i,j)$.  One then builds $P$ and $Q$ recursively by adding each such $i$ to $P$ in a particular way such that the result is always an SSYT, with the entries $j$ of the biletters becoming the entries of $Q$. The procedure for adding each $i$ to $P$ is known as \emph{RSK insertion}.      

To generalize RSK insertion to skew shapes, two insertion procedures are defined in \cite{SS90}.  Starting with an SSYT $T$ of shape $\lambda/\mu$, the first type of insertion, \emph{external insertion}, works just like RSK insertion.  To describe \emph{internal insertion}, we first say that a cell $(a,b)$ of $\lambda/\mu$ is an \emph{inner corner} if $(a,b-1), (a-1,b) \not\in \lambda/\mu$.  Internal insertion works by removing the entry in such a cell $(a,b)$ and inserting it into row $a+1$ using the usual RSK insertion procedure.  Note that internal insertion, unlike external insertion, does not increase the number of entries of $T$.  

These insertions are used to prove a number of bijections, including a bijection from tuples $(A, T, U) $ to $(P,Q)$, where $A=(a_{ij})$ is a matrix over $\mathbb{N}$ and $T, U, P, Q$ are SSYT of shape $\alpha/\mu$, $\beta/\mu$, $\lambda/\beta$ and $\lambda/\alpha$ respectively.  Here, $\alpha$ and $\beta$ are fixed partitions, while $\lambda$ depends on the choice of $(A, T, U)$.

Just like the Cauchy identity follows from the classical RSK algorithm, the following generalization follows from the bijection just described.  For fixed partitions $\alpha$ and $\beta$,
\begin{equation}\label{equ:gencauchy}
\sum_\lambda s_{\lambda/\beta}(X) s_{\lambda/\alpha}(Y) = \sum_\mu s_{\alpha/\mu}(X) s_{\beta/\mu}(Y) \prod_{i,j} (1-x_i y_j)^{-1}.
\end{equation}
Independent proofs of this identity using symmetric function techniques have been given by Lascoux, Macdonald, Towber, and Zelevinsky \cite{Mak85}, \cite[Example~I.5.26]{Mac95}.
The resulting analogue of \eqref{equ:flambda} is as follows, where $n$ and $m$ are fixed integers and $\alpha$ and $\beta$ are again fixed partitions:
\[
\sum_{\genfrac{}{}{0pt}{}{\lambda/\beta \vdash n}{\lambda/\alpha \vdash m}} f_{\lambda/\beta} f_{\lambda/\alpha} = \sum_{k \geq 0} \binom{n}{k} \binom{m}{k} k! \sum_{\genfrac{}{}{0pt}{}{\alpha/\mu \vdash n-k}{\beta/\mu \vdash m-k}} f_{\alpha/\mu} f_{\beta/\mu} .
\]
Notice that letting $n=m$ and $\alpha$ and $\beta$ be empty yields \eqref{equ:flambda}.  Equation \eqref{equ:gencauchy} is one of eight identities which make up \cite[Exer.\ 7.27]{ec2}, all of which appear in \cite{SS90} as consequences of their various bijections.

\subsection{Recent applications}

Our presentation of some applications of the insertion procedures of \cite{SS90} begins with a conjecture of Richard from \cite{Sta05}.  Let the \emph{sign} of an SYT (in English notation) be the sign of the permutation obtained by reading the rows from left to right, starting with the top row.  The \emph{sign imbalance} $I_\lambda$ of a partition $\lambda$ is the sum of the signs of all SYT of shape $\lambda$, and significant attention has been given to the question of determining the sign imbalance of partition shapes and characterizing those with $I_\lambda=0$.  Richard conjectured \cite[Conjecture~3.3(a)]{Sta05} that the sum of $I_\lambda$ over all shapes $\lambda$ with $n$ cells is $2^{\lfloor n/2 \rfloor}$.  This conjecture was subsequently proved independently by Lam, Reifegerste and Sj\"ostrand \cite{Lam04, Rei04, Sjo05}.  The technique used by Reifegerste and Sj\"ostrand is to establish the relationship between the sign of a permutation $\pi$ and 
the sign of the image $(P,Q)$ of the permutation matrix for $\pi$ under the classical RSK algorithm.  

In \cite{Sjo07}, Sj\"ostrand asks about the sign imbalance of \emph{skew} shapes $\lambda/\mu$.  Let $A$ denote the permutation matrix of a permutation $\pi$, and suppose $(A,T,U)$ is mapped to $(P,Q)$ under the bijection of \cite{SS90}.\footnote{Sj\"ostrand worked with the more general notion of a \emph{partial permutation} $\pi$, but, for simplicity, we restrict our attention to the case when $\pi$ is a permutation.}  Sj\"ostrand's main theorem establishes a remarkably simple relationship between the signs of $\pi$, $T$ and $U$ and those of $P$ and $Q$.  This relationship is then used to establish skew analogues of identities from the partition case.  It will come as no surprise to the reader that the key relationship is established using the external and internal insertions of \cite{SS90}.

In \cite{AM11}, Assaf and the second author generalize the Pieri rule to skew shapes by giving an expansion of $s_{\lambda/\mu}s_{(n)}$ as a signed sum of skew Schur functions.  The proof is combinatorial, and the insertion procedures of \cite{SS90} are exactly what is needed to establish the crucial sign-reversing involution.  The same applies to the proofs in \cite{Kon12}, where Konvalinka presents a simpler involution that proves a dual version of the skew Pieri rule.  In addition, he combines his involution with the one in \cite{AM11} to prove a ``skew quantum Murnaghan--Nakayama rule,'' which simultaneously generalizes the classical Murnaghan--Nakayama rule, the skew Pieri rule and its dual, and other related results.  

Three other recent articles follow this same thread.  Assaf and the second author include a conjecture of a ``skew Littlewood--Richardson rule,'' i.e., an expansion of $s_{\lambda/\mu}s_{\sigma/\tau}$ as a signed sum of skew Schur functions.  This conjecture was proved by Lam, Lauve and Sottile \cite{LLS11}, along with several other results of a similar flavor.  Konvalinka and Lauve \cite{KL13} provide skew Pieri rules for Hall--Littlewood functions, thereby introducing a parameter $t$ into the story.  Specifically, they give an expression for the product of the skew Hall--Littlewood polynomial $P_{\lambda/\mu}$ times $h_r$ as a signed sum of skew Hall--Littlewood polynomials, and do the same with $e_r$ or $q_r := (1-t)P_r$ in place of $h_r$.  Finally, in \cite{War13}, Warnaar shows that $q$-analogues of these three results from \cite{KL13} can be derived from a $q$-binomial theorem for Macdonald polynomials of Lascoux and himself \cite{LW11}.

Beyond \cite{SS90}, Richard played an additional role in initiating the results of the previous two paragraphs.  Assaf and the second author were both at MIT in the spring of 2009 when they stumbled upon the possibility of a simple expansion for $s_{\lambda/\mu}s_{(1)}$.  The obvious thing to do in such a situation is to ask Richard if it is already known.  Richard was surprised by the expansion but that same day provided an algebraic proof of this $n=1$ case.  He provided encouragement as the conjecture for $s_{\lambda/\mu}s_{(n)}$ in the general $n$ case was formulated and then given a combinatorial proof.

\section{Flag symmetry of posets}\label{sec:flag}

In Section~\ref{sec:chromatic}, we asked for naturally-arising symmetric functions and looked at a topic that lies at the intersection of graph theory and symmetric function theory.  This section has a similar flavor, but now the symmetric function of interest is defined in terms of chains in posets. 
Specifically, we begin with the following quasisymmetric function introduced by Richard Ehrenborg in \cite{Ehr96}.  For a finite ranked poset $P$ with $\hat{0}$ and $\hat{1}$ and rank function $\rho$\,, define a formal power series in the variables $X=(x_1, x_2, \ldots)$ by
\[
\Ehr_P(X) = \sum_{\hat{0}=t_0 \leq t_1 \leq \cdots \leq t_{k-1} < t_k =\hat{1}} x_1^{\rho(t_0,t_1)}x_2^{\rho(t_1,t_2)}\cdots x_k^{\rho(t_{k-1},t_k)},
\]
where $\rho(t_{i-1}, t_i)$ denotes $\rho(t_i)-\rho(t_{i-1})$, and where 
 the sum is over all multichains from $\hat{0}$ to $\hat{1}$ such that $\hat{1}$ occurs exactly once.  (This last requirement ensures that the coefficients of $\Ehr_P(X)$ are finite.)  

\subsection{Flag symmetry}

Ehrenborg asked for which posets $P$ is $\Ehr_P(X)$ symmetric, and Richard (Stanley!) has three papers \cite{Sta96, Sta97, SS99} that address this question, the third of which is joint work with Rodica Simion.
In \cite{Sta96}, Richard terms $P$ ``flag symmetric" if $\Ehr_P(X)$ is symmetric; the reason for the terminology is that $\Ehr_P(X)$ encodes the same information as the flag $f$-vector and flag $h$-vector of $P$.  Indeed, if $P$ has rank $n$, let $S = \{m_1, \ldots, m_j\}_< $ denote a subset of $[n-1]$ satisfying $m_1 < \cdots < m_j$.  
Then the flag $f$-vector $\alpha_P(S)$ is defined as the number of chains $\hat{0} < t_1 < \cdots < t_j =\hat{1}$ such that $\rho(t_i) = m_i$ for all $i$.  The flag $h$-vector $\beta_P(S)$ is defined by $\alpha_P(S) = \sum_{T \subseteq S} \beta_P(T)$.
Richard gives the following alternative expressions for $\Ehr_P(X)$:
\begin{eqnarray*}
\Ehr_P(X) &  = & \sum_{\genfrac{}{}{0pt}{}{S = \{m_1, \ldots, m_j\}_<}{S \subseteq [n-1]}} \sum_{1 \leq i_1 < \cdots < i_{j+1}} \alpha_P(S)\,
x_{i_1}^{m_1} x_{i_2}^{m_2-m_1} \cdots x_{i_{j+1}}^{n-m_j} \\
& = & \sum_{S \subseteq [n-1]} \beta_P(S) \Fund_{S}(X).
\end{eqnarray*}

Richard's main tool for showing rank symmetry is the notion of local rank symmetry, which Richard proves is sufficient for flag symmetry.  A poset is said to be \emph{locally rank symmetric} if all of its intervals are rank symmetric, meaning that the cardinalities of the rank levels read from bottom to top form a palindromic sequence.  Besides a lot of other results, Richard gives the following examples of flag symmetric posets:
\begin{itemize}
\item a finite distributive lattice is flag symmetric if and only if it is a product of chains, in which case $\Ehr_P = h_\nu$, the complete homogeneous symmetric function indexed by the chain lengths \cite{Sta96};
\item binomial posets, observed to be flag symmetric in \cite{Ehr96};
\item face lattices of simplices, polygons, or of three-dimensional polytopes with equal numbers of vertices as facets, as well as the product of any of these types of face lattices \cite{Sta96};
\item the lattice of non-crossing partitions \cite{Sta97}; 
\item the poset of shuffles as introduced by Curtis Greene \cite{Gre88}, shown to be locally rank symmetric in \cite{SS99}.
\end{itemize}

\subsection{Representation theory of the symmetric group}

A second thread of results that we choose to highlight among all those in \cite{Sta96, Sta97, SS99} is the elegant ways in which $\Ehr_P(X)$ relates to the representation theory of the symmetric group.  Although the requirements for these relationships to hold are very special, Richard's work gives a number of instances where everything works out beautifully.  When $\Ehr_P(X)$ is symmetric, it can be expanded in terms of the basis of Schur functions:
\[
\Ehr_P(X) = \sum_\lambda a_\lambda s_\lambda(X).
\]
Recall that if all the $a_\lambda$ are nonnegative, then $\Ehr_P(X)$ is said to be \emph{Schur positive}.  As a result, $\Ehr_P(X)$ must equal the Frobenius characteristic $\ch(\psi)$ of some character $\psi$ of the symmetric group (see \cite[\S 7.18]{ec2}).  Richard then seeks a natural action of $S_n$ on the maximal chains of $P$ that would give rise to $\psi$.  Such an action should be \emph{local}, meaning that when the transposition $(i,i+1)$ acts on a maximal chain $m$, the result should be a linear combination of maximal chains that differ from $m$ only at rank $i$.  Putting these ideas together gives the definition of what Richard calls a ``good'' action: a local action on the maximal chains of $P$ with the property that $\Ehr_P(X) = \ch(\psi)$ or potentially $\omega(\Ehr_P(X)) = \ch(\psi)$, where $\omega$ is the usual involution on symmetric functions.  

The three papers we have been discussing exhibit good $S_n$ actions in the following cases.
\begin{itemize}
\item A product of chains \cite{Sta96}.  In this case, the adjacent transposition $(i,i+1)$ of $S_n$ sends a maximal chain $m$ to the unique chain $m'$ that differs from $m$ only at rank $i$, while $m$ remains fixed if no such $m'$ exists.  Conversely, it follows from a result of Grabiner \cite{Gra99} that if $P$ has a good action of this form and additionally is a Cohen--Macaulay poset (defined in \cite{BGS82} or in \S3.8 of \cite{ec1,ec1e2}), then $P$ is a product of chains.  
\item The lattice of non-crossing partitions \cite{Sta97}.  In this case, $\omega(\Ehr_P(X))$ equals Haiman's parking function symmetric function, defined as the Frobenius characteristic arising from the $S_n$ action on parking functions that permutes coordinates \cite{Hai94}.  Richard
 thus establishes a remarkable connection between non-crossing partitions and parking functions.  This connection has a number of nice consequences, e.g., the lattice of non-crossing partitions of $\{1,\ldots,n+1\}$ can be given an edge-labeling so that the maximal chains are labeled by the parking functions of length $n$, each occurring exactly once.  
\item The poset of shuffles \cite{SS99}.  In fact, \cite{SS99} shows the more general result that a good action results anytime a poset has a chain labeling with certain properties. 
\end{itemize}

\section{Conclusion}\label{sec:conclusion}

Hopefully, this paper inspires further interest in symmetric and
quasisymmetric functions by highlighting some the common threads and
interactions among the results.  We have just scratched the surface of
all the amazing mathematical contributions, which will be recognized
for generations to come, by YLNTAES RCHRADI.

\section*{Acknowledgments}\label{sec:ak}

On a more personal note, we are forever grateful to Richard Stanley
for mentoring both of us as graduate students and throughout our
careers.\footnote{The second author was officially Richard's
  Ph.D. student.  The first author was a visiting student at MIT for
  1.5 years of graduate school.  She considers herself to be a
  Ph.D. foster student of Richard's.  He was also her NSF postdoc
  mentor.}  We appreciate the clarity of his ideas, his enthusiasm for
research, and his enormous wealth of mathematical knowledge.  We also
appreciate his fun sense of humor, hence we added a few extra jokes
which we hope the readers enjoyed.

Many people read and commented on earlier drafts of this paper and we
are very grateful for their insights.  In particular, we would like to
thank Matja\v z Konvalinka, Ezra Miller, Richard Stanley, Josh Swanson, Jair Taylor, and Vasu Tewari for many helpful
suggestions and comments.

\bibliographystyle{alpha}
\bibliography{rstan_sym}

\end{document}